%% file: main.tex
\documentclass[preprint,11pt]{elsarticle}
\usepackage[T1]{fontenc}
\usepackage{lmodern,mathtools,amssymb,amsthm,booktabs}
\usepackage[hidelinks]{hyperref}
\usepackage[capitalise,noabbrev]{cleveref}
\usepackage{microtype}
\input{metadata}
\newtheorem{theorem}{Theorem}[section]
\newtheorem{lemma}[theorem]{Lemma}
\newtheorem{prop}[theorem]{Proposition}
\crefname{prop}{Proposition}{Propositions}
\newcommand{\R}{\mathbb R}
\newcommand{\E}{\mathbb E}
\newcommand{\tr}{\operatorname{tr}}
\newcommand{\diag}{\operatorname{Diag}}
\newcommand{\norm}[1]{\lVert #1\rVert}
\newcommand{\ip}[2]{\langle #1,#2\rangle}
\newcommand{\T}{^{\mathsf T}}
\begin{document}
\begin{frontmatter}
\title{\PaperTitle}
\author[btbu]{\ZhaoName}
\author[ustc]{\LiName}
\address[btbu]{\ZhaoAff, \ZhaoPostal}
\address[ustc]{\LiAff, \LiPostal}
\begin{abstract}
Consider $n$ unit vectors in a Euclidean space of dimension $d\ge1$. We prove that some signed sum has length at most $\sqrt{d-1}$ when $n$ and $d$ have opposite parity, as conjectured by Hollom, Portier and Souza. For $n\ge d$, both this bound and the classical bound $\sqrt d$ for equal parity are sharp. When $n<d$, the sharp bound is $\sqrt n$ in either case. For $d\ge3$, a counterexample would give an odd-dimensional central section of a cube. We rule out the possible dimensions by face incidences and, in the last case, by a Gram matrix rank bound.
\end{abstract}
\end{frontmatter}

\section{Introduction}
The Euclidean discrepancy of a family of vectors is the length of its shortest signed sum. For $V=[v_1,\ldots,v_n]\in\R^{d\times n}$, write
$D(V)=\min_{\varepsilon\in\{\pm1\}^n}\norm{V\varepsilon}^2$
for its square. With independent uniform signs, the reverse Littlewood--Offord problem asks how small the probability of hitting a closed ball about the origin can be \cite{hps,hs}. For a fixed family and radius $r\ge0$, this probability is positive exactly when $D(V)\le r^2$. At radius $\sqrt d$, Beck's probability bound \cite{beck} gives $D(V)\le d$ for every family of unit vectors (see also \cite[Theorem~1.1]{hs}).

When $n\not\equiv d\pmod2$ and $d\ge3$, Hollom and Sorkin proved $D(V)\le d-\epsilon(d)$ with $\epsilon(d)=2^{-200}d^{-112}$ \cite[Theorem~1.2]{hs}. Hollom, Portier and Souza conjectured that $D(V)\le d-1$ \cite{hps} (see \cite[Conjecture~1.5]{hs}). We prove this conjecture in every dimension. Hollom and Sorkin handle families clustered near orthogonal directions and use a stability estimate when two vectors are neither nearly parallel nor nearly orthogonal \cite[Sections~4 and~5]{hs}. If every signed sum has squared norm greater than $d-1$, we compress signed sums on even supports. The resulting cube section has odd dimension, and its face incidences and a rank bound give a contradiction.

\begin{theorem}\label{th:main}
Let $d\ge1$, $n\ge0$ and $n\not\equiv d\pmod2$. If $\norm{v_i}=1$ for every $i$, then some $\varepsilon\in\{\pm1\}^n$ satisfies
\[
 \Big\lVert\sum_{i=1}^n\varepsilon_i v_i\Big\rVert^2\le d-1 .
\]
\end{theorem}

For fixed $n$ and $d$, the largest value of $D(V)$ over all families of unit vectors is
\[
 \max_V D(V)=
 \begin{cases}
 n   & n<d,\\
 d   & n\ge d,\quad n\equiv d\pmod2,\\
 d-1 & n\ge d,\quad n\not\equiv d\pmod2 .
 \end{cases}
\]
For $n<d$, independent fair signs give $D(V)\le\E\norm{V\varepsilon}^2=n$, with equality for $n$ orthonormal vectors. For $n\ge d$, take $e_1,\ldots,e_{d-1}$ and $n-d+1$ copies of $e_d$. The first $d-1$ coordinates have magnitude one under every signing. The minimum magnitude of the last coordinate is one in equal parity and zero in opposite parity. This gives equality in Beck's bound and in Theorem~\ref{th:main}, respectively.

Thus $\epsilon(d)=1$ is optimal for every $d\ge3$. For fixed $d\ge3$ and $n$ of opposite parity, Hollom and Sorkin give examples in which the probability that a random signed sum lies in the closed ball of radius $\sqrt{d-1}$ is at most $2^{-n/2+O_d(1)}$ \cite[Section~6]{hs}. Can the probability decay at a faster exponential rate for some choice of vectors?

In dimension one an even number of unit vectors can be signed to cancel in pairs. Swanepoel proved the planar bound for odd families of unit vectors in any normed plane \cite[arXiv v1, Theorem~2]{swanepoel}. Embedding the plane in three dimensions gives another proof of the Euclidean planar case.

Vertex rounding gives the bound $d$ in any $d$-dimensional normed space. At most $d$ coefficients are free at a vertex of the cube section. Rounding one such coefficient to its nearest sign changes the sum by at most one in norm. Swanepoel suggested that $d-1$ might suffice for $n$ unit vectors when $n\not\equiv d\pmod2$ \cite[arXiv v1, Section~4]{swanepoel}. In $\ell_2^d$ with $d\ge3$, Beck's bound suffices since $\sqrt d\le d-1$. Theorem~\ref{th:main} gives the sharp bound $\sqrt{d-1}$ for opposite parity.

In $\ell_1^d$ the following family attains $d-1$. For $n=d-1+2k$ with $k\ge0$, take $e_1,\ldots,e_{d-1}$ and $2k$ copies of $e_d$. The first $d-1$ coordinates have magnitude one, and the copies of $e_d$ can cancel in pairs. The shortest signed sum therefore has length $d-1$ in $\ell_1^d$ and $\sqrt{d-1}$ in $\ell_2^d$. In \emph{Parity Balancing of $\ell_1$ Unit Vectors} (submitted to arXiv), we prove the matching upper bound in $\ell_1^d$ for every $d\ge2$.

\subsection{Proof strategy}
Suppose $n$ and $d$ have opposite parity. We seek a signed sum of length at most $\sqrt{d-1}$. First allow coefficients $x_i$ in $[-1,1]$ and look for exact cancellation, $\sum_i x_i v_i=0$. The solutions form a section of the cube. At each vertex, at most $d$ coefficients lie strictly between $-1$ and $1$. All other coefficients are signs. If we round the remaining coefficients to signs $\varepsilon_i$, then $V\varepsilon=V(\varepsilon-x)$, so the signed sum is exactly the rounding error. With at most $d-1$ coefficients to round, \eqref{eq:round-cost} gives the required bound.

Suppose $d\ge3$ and every signing fails. We shorten the list by replacing disjoint groups of even size by short signed sums. Either sign on a group sum lifts to a signing of that group. Thus every signing of the shorter list still fails. The rounding bound on this list also constrains the vectors left outside the groups. Their cube section has positive odd dimension $r$. At each vertex exactly $r$ coefficients are already signs, and the sum of the corresponding vectors with those signs has length less than one.

This rules out $r=1$, since a single unit vector has length one. If $r\ge5$, two opposite coordinate walls would cut out facets of the same three-dimensional face. Lemma~\ref{lem:axis} shows that its facets must lie on distinct axes. If $r=3$, the triples of fixed signs involve only six coordinates. Their relative signs determine a five-cycle. Combining columns according to this pattern gives a matrix $Y$ with $d+4$ columns in $\R^d$. Its kernel therefore has dimension at least four. The Gram matrix estimate for the five-cycle bounds the same dimension by three.

\subsection{Rounding and the terminal cube section}
Call $\sum_{i\in I}\sigma_i v_i$, with $\sigma_i=\pm1$, a \emph{word} on $I$. Its parity is that of $|I|$. For the cube-section argument, assume $d\ge3$ and $n\not\equiv d\pmod2$. Suppose for a contradiction that
\begin{equation}\label{eq:H}
  D(V)>L,\qquad L=d-1 .
\end{equation}
Compress words on disjoint supports into single columns. Every signing of the resulting matrix $W$ lifts to a signing of $V$, so $D(W)>L$. Set $K(W)=\{x\in[-1,1]^m:Wx=0\}$ and $F(x)=\{i:|x_i|<1\}$. The coordinates in $F(x)$ are \emph{free} and the others are \emph{saturated}.

At a vertex, fixing the saturated walls leaves no kernel direction on the free coordinates. Under the hypotheses below exactly $m-d$ walls are saturated.

The free columns form a basis and hence determine the number of saturated walls. The two bounds in \eqref{eq:def} control the saturated word and rounding on small coordinate sets.

\begin{lemma}[Two rounding rules]\label{lem:round}
Suppose $W=[w_1,\ldots,w_m]$, $\norm{w_i}\le1$, and $D(W)>L$. At a vertex $x$ of $K(W)$ write $F=F(x)$ and $z=\sum_{i\notin F}x_iw_i$. Then the columns of $W_F$ form a basis of $\R^d$, the polytope $K(W)$ is simple, and
\begin{equation}\label{eq:def}
 \norm z^2+\sum_{i\in F}(1-\norm{w_i}^2)<1,
 \qquad \sum_{i\in I}|x_i|<1\quad(I\subseteq F,\ |I|\le5).
\end{equation}
\end{lemma}
\begin{proof}
At a vertex, $W_F$ has independent columns. Fix the saturated signs and switch column signs so $x_F\ge0$ coordinatewise. On $I\subseteq F$ put $t_i=1-x_i$. If $\sum_I t_i\le4$, give $I$ positive signs and flip at most one, choosing $i$ with probability $t_i/4$. Then $h=(\varepsilon-x)_I$ satisfies $\E hh\T=\operatorname{diag}(t_i)_{i\in I}$. For each remaining free sign, apply \cite[Lemma~3.2, (3.1)]{hs}. The resulting choice increases the squared error by at most $(1-x_i^2)\norm{w_i}^2$. Thus
\begin{equation}\label{eq:round-cost}
 D(W)\le\E\norm{W(\varepsilon-x)}^2
 \le |F|-\sum_{i\in I}x_i-\sum_{i\in F\setminus I}x_i^2 .
\end{equation}
With $I=\varnothing$, $d-1<D(W)\le |F|-\norm{x_F}^2\le d$, hence $|F|=d$ and $\norm{x_F}^2<1$. The complementary coordinates parametrize the kernel, so $K(W)$ is simple. If $|I|\le5$ and $\sum_I x_i\ge1$, then $\sum_I t_i\le4$ and \eqref{eq:round-cost} gives $D(W)\le d-1$.

For $G=Z\T Z$, let $\diag G$ be its diagonal part. The second and fourth moments of independent fair signs give
\begin{equation}\label{eq:moment}
 \E[(\varepsilon\T G\varepsilon-t)\varepsilon\varepsilon\T]
 =2(G-\diag G)+(\tr G-t)I .
\end{equation}
If $D(Z)>t$, the left-hand side is positive definite. Compress the saturated columns into $z$. Take $Z=(z,W_F)$, $t=L$, and put $e_i=1-\norm{w_i}^2$ and $u=\norm{x_F}^2<1$. Evaluating \eqref{eq:moment} on the kernel vector $(1,x_F)$ gives
\[
 (1-u)\left(1-\norm z^2-\sum_{i\in F}e_i\right)
 >2\sum_{i\in F}e_i(u-x_i^2)\ge0 .
\]
\end{proof}

Lemma~\ref{lem:round} gives $\operatorname{rank}V=d\le n$. Thus $n=d+k$ with $k$ positive and odd. Start with $(r,S)=(k,0)$.
Each removed support has even size at most $r-1$, so the excess $r$ stays positive and odd.
Among nonempty even words on at most $r-1$ remaining columns, take a shortest one $E$. Stop if there is no such word or if $S+\norm E^2>1$. Otherwise remove its support and update
\[
 (r,S)\longleftarrow(r-|\operatorname{supp}E|,\ S+\norm E^2).
\]
Write $E_1,\ldots,E_j$ for the selected words and $A$ for the remaining columns.

A shortest word minimizes the increment of $S$. If its block were free at a vertex, the saturated columns would yield a shorter competitor.

\begin{prop}[The terminal section]\label{prop:kernel}
We have $S<1$. Put $R^2=1-S$. Then $K=K(A)$ is a simple centrally symmetric $r$-polytope. Each vertex saturates exactly $r$ axes, and its saturated word $F$ satisfies
\begin{equation}\label{eq:corners}
 \norm F^2<R^2,\qquad \Big\lVert F+\sum_b\theta_b E_b\Big\rVert^2<1\quad(\theta_b=\pm1).
\end{equation}
The five-coordinate bound in \eqref{eq:def} holds at each vertex. Vertex words on the same support agree up to a global sign. Every nonempty even word on fewer than $r$ remaining columns has norm greater than $R$.
\end{prop}
\begin{proof}
Write $E=[E_1,\ldots,E_j]$ and $\widehat K=K([A,E])$. Lemma~\ref{lem:round} applies since each block has squared norm at most $S\le1$. Two free blocks would have total deficit at least $2-S\ge1$, contrary to \eqref{eq:def}. If a block is saturated at every vertex, release its wall along an incident edge. By simplicity, all other active walls remain fixed along the edge. At the other endpoint the block is again saturated, with the opposite sign. Thus one can flip any such block while preserving all other saturated signs and the free index set.

Suppose $i$ is the least index for which $E_i$ is free at some vertex, and fix such a vertex. Its free columns are $E_i$ and $d-1$ columns of $A$. Omit the earlier blocks and expand the remaining saturated columns into a word $C$ on the original vectors. Just before step $i$ there were $d+t$ vectors left, where $t$ was the excess. Of these, only the support of $E_i$ and the $d-1$ free singletons are missing from $C$. Hence
\[
 |\operatorname{supp} C|=t+1-|\operatorname{supp}E_i|\in\{2,4,\ldots,t-1\}.
\]
Thus $C$ was eligible when we chose $E_i$. Minimality gives $\norm{E_i}\le\norm C$.

All blocks $E_b$ with $b<i$ are saturated everywhere. Flipping them realizes every saturated word $C+\sum_{b<i}\theta_bE_b$ while $E_i$ remains free. This free block contributes $1-\norm{E_i}^2$ to the deficit in \eqref{eq:def}, so $\norm{C+\sum_{b<i}\theta_bE_b}^2<\norm{E_i}^2$ for every choice of signs. For independent fair $\theta_b$ the cross terms have mean zero. Averaging gives
\[
 \norm{E_i}^2\le\norm C^2\le\norm C^2+\sum_{b<i}\norm{E_b}^2
 =\E_\theta\Big\lVert C+\sum_{b<i}\theta_b E_b\Big\rVert^2<\norm{E_i}^2 .
\]
Hence every block is saturated at every vertex.

Fix a vertex $x$ of $K$ and let $P$ be the face of $\widehat K$ cut out by its saturated walls. For every vertex $y$ of $P$, all blocks are saturated and $F(y)\subseteq F(x)$. Since $|F(y)|=d$ and $A_{F(x)}$ has independent columns, $F(y)=F(x)$ and these columns form a basis. Thus $K$ is simple of dimension $r$. The face contains $(x,0)$, which inherits the five-coordinate bound by convexity. Independent block flips and \eqref{eq:def} give $\norm{F+E\theta}^2<1$ for every $\theta$. Averaging gives $\norm F^2+S<1$, proving \eqref{eq:corners}.

The stopping rule gives the lower bound for even words. If vertex words $F,F'$ on the same support satisfy $F'\ne\pm F$, choose the sign of $F'$ so the number of disagreements is even ($r$ is odd). Then $(F-F')/2$ is a nonempty even word on fewer than $r$ columns with norm less than $R$, a contradiction.
\end{proof}

\section{Parallel edges and low-dimensional faces}
We restrict the walls of a three-face, then pass to higher dimensions through vertex links. Two opposite facets of the section would descend to opposite walls on a single three-face.
A wall is an equation $x_i=\pm1$, and its axis is $i$. An $s$-face of $K$ fixes $r-s$ walls. The vertices and maximal simplices of its dual boundary correspond to its facets and vertices, respectively.

At a vertex $x$, let $s_i(x)=x_i\mathbf1_{|x_i|=1}$. If $I$ contains at most five free coordinates, then $\norm{(x-s(x))_I}_1<1$. In particular, suppose a face $P$ contains vertices $v_1,\ldots,v_k$ and the axes $i_1,\ldots,i_k$ are distinct, with $k\le6$. If $v_a$ saturates only $i_a$ among these axes, then
\begin{equation}\label{eq:private}
 k\le\dim P+1,
\end{equation}
since the evaluation matrix $((v_a)_{i_b})$ is strictly diagonally dominant and has rank at most $\dim P+1$.

\begin{lemma}[Common directions and facet axes]\label{lem:sign}
Write $B(P)$ for the facet axes of $P$ and $\operatorname{dir}P=\operatorname{span}(P-P)$. For one- or two-dimensional faces $P,Q$ of $K$,
\[
 \operatorname{dir}P\cap\operatorname{dir}Q\ne\{0\}
 \quad\Longrightarrow\quad B(P)\cap B(Q)\ne\varnothing .
\]
If $r\ge3$, two vertex words sharing at least $r-1$ axes have the same signs on the common axes up to a global sign.
\end{lemma}
\begin{proof}
Suppose the axes are disjoint and $0\ne w\in\operatorname{dir}P\cap\operatorname{dir}Q$. In each $R=P,Q$ take a chord $y-v=\lambda_R w\ne0$ from a vertex to a point on a facet $e$ that does not contain the vertex. Take a parallel edge if there is one. Otherwise take a maximal chord through a vertex whose transverse projection lies in the interior of the projected interval. Let $I_R$ consist of the axes of $e$ and the facets at $v$, and set $I=I_P\sqcup I_Q$.

Write $y$ as a convex combination of vertices of $e$ and use the same coefficients for $s(y)$. Put $u=(s(y)-s(v))_I$. Outside $B(R)$ a coordinate is fixed on $R$ or free at every vertex of $R$, so $\operatorname{supp}u\subseteq I_R$. On a wall of $v$ with sign $\sigma$, the change in the saturated entry is $0,-\sigma$ or $-2\sigma$. On the axis of $e$ the change has the sign of its wall. Thus interpolation gives no cancellation on $I_R$. Each vertex of $e$ gains the wall of $e$ and loses at least one wall through $v$, since the facets of $R$ through $v$ intersect only at $v$. If the gained and lost walls have the same axis, the change there has magnitude two. Hence $\norm u_1\ge2$.

At $v$ and at each vertex used in the interpolation, at least one of the at most six axes in $I$ is saturated. By \eqref{eq:def} the free entries on $I$ have $\ell_1$ norm less than one at each of these vertices. The same bound holds after interpolation. Taking the difference of the endpoint errors gives
\[
 \norm{\lambda_R w_I-u}_1<2\le\norm u_1 .
\]
On $I\setminus I_R$ we have $u=0$. The reverse triangle inequality on $I_R$ now gives $\norm{w_{I_R}}_1>\norm{w_{I\setminus I_R}}_1$. Taking $R=P,Q$ gives opposite strict inequalities.

Equal supports are covered by Proposition~\ref{prop:kernel}. For $r-1$ common axes $J$, fixing their signs at the two vertices gives parallel edges with direction $\{h\in\ker A:h_J=0\}$. A common facet axis gives endpoints of equal support, so the signs on $J$ agree up to negation.
\end{proof}

\begin{lemma}[Axes of low-dimensional faces]\label{lem:axis}
If $r\ge3$, a two-face has at most five edges, with distinct axes. If $r\ge5$, the facets of a three-face have distinct axes, and any two disjoint facets have a common neighboring facet.
\end{lemma}
\begin{proof}
Vertices of a proper face agree on its fixed walls. Thus the global sign in Lemma~\ref{lem:sign} must be positive whenever two such vertices share at least $r-1$ axes. Opposite edges on one axis of a two-face would give endpoints sharing at least $(r-2)+1$ axes with a sign disagreement. If the polygon had at least six edges, number them in cyclic order. The initial vertices of edges $1,4$ and the terminal vertices of edges $2,5$ each lie on just one of these four edges. Their four distinct axes would give $4\le3$ in \eqref{eq:private}.

For disjoint facets $P,Q$ of a three-face, the two direction planes lie in the same three-dimensional space and have a nonzero vector in common. Take $j\in B(P)\cap B(Q)$ by Lemma~\ref{lem:sign}. If $P,Q$ were opposite on one axis, their $j$-edges would have endpoints sharing at least $(r-3)+2$ axes with a sign disagreement. Thus facet axes are distinct, and the facet on $j$ meets both.
\end{proof}

Assume $r\ge5$ for the rest of the section. Let $\Delta$ be the dual boundary of a three-face and $N[v]=N(v)\cup\{v\}$ the closed neighborhood. A facial triangle $T$ is invisible from $v$ if $T\cap N[v]=\varnothing$.
\begin{lemma}[Invisible triangles]\label{lem:local}
If a facial triangle is invisible from $v$, then $\Delta$ has seven vertices, $v$ has degree three, and this triangle is unique. There are at most four degree-three vertices, with equality only for the boundary of a tetrahedron.
\end{lemma}
\begin{proof}
Put $G=\Delta^{(1)}$, $n=|V(G)|$ and $W=V(G)\setminus N[v]$. For $w\in W$, take a common neighbor $x$ of $v,w$ by Lemma~\ref{lem:axis}. In the link cycle of $x$ (length at most five), $v,w$ have distance two. The link path $v,y,w$ gives a facial triangle $xyw$ with $x,y\in N(v)$.

Suppose $W$ contains four vertices $w_1,\ldots,w_4$. Among the five labels $v,w_1,\ldots,w_4$, the triangle constructed for $w_i$ contains only $w_i$. A triangle through $v$ contains only $v$ from this list. The labels correspond to distinct facet axes by Lemma~\ref{lem:axis}. Each of the five dual triangles thus gives a vertex of the primal three-face that saturates only the corresponding axis among the five. Now \eqref{eq:private} with $k=5$ gives $5\le3+1$, a contradiction. Thus $|W|\le3$ and $\deg v\ge n-4$ for all $v$. Summing the degrees and using Euler's formula gives $n(n-4)\le6n-12$, so $n\le8$. Also $\sum_v(6-\deg v)=12$ bounds the number of cubic vertices by four, with equality only for the tetrahedron.

Orient the facet walls positively. For colors $c_i=\pm1$, the three saturated coordinates of the primal vertex $x_F$ contribute $\sum_{i\in F}c_i$ to $\sum_i c_i(x_F)_i$. The first sum is an odd integer. The other $n-3\le5$ coordinates contribute less than one in absolute value by \eqref{eq:def}, so the functional is positive precisely at the vertices whose triangles have a positive majority. The positive vertices in the primal graph have increasing paths to the maximizing face, whose graph is connected \cite{kalai}. These paths stay positive, so their dual triangles are edge-connected.

Suppose four vertices induce $2K_2$ and color just those four positive. A triangle has positive majority exactly when it contains one of the two positive edges. Each edge lies in a facial triangle. Triangles containing different positive edges cannot share an edge, since no endpoint of one positive edge is adjacent to an endpoint of the other. This contradicts connectivity.

An invisible triangle equals $W$. Each of the $k=\deg v$ neighbors is adjacent to at least two of its vertices, since an uncovered edge would form an induced $2K_2$ with the edge to $v$. Counting edges incident to $v$, edges of its link, edges in $W$, and edges from the link to the triangle gives $4k+3\le3(k+4)-6$. Thus $k=3$ and $n=7$.
\end{proof}

\begin{prop}[Excluding higher dimensions]\label{prop:high}
The case $r\ge5$ is impossible.
\end{prop}
\begin{proof}
In the dual boundary $Q$ of a face, call a maximal simplex invisible from $v$ if it avoids $N[v]$. Suppose some maximal simplex is invisible and let $B$ be the boundary over $\mathbb F_2$ of the sum of all invisible maximal simplices. This sum omits every maximal simplex through $v$. By ridge connectivity, some ridge $T$ separates an invisible maximal simplex from a visible one. That ridge occurs in the boundary, so $B\ne0$. For every such $T$, the visible simplex is $T\cup\{x\}$ with $x\in N[v]$. Since $T$ has no vertex adjacent to $v$, we have $x\ne v$. Thus $x\in N(v)$ and $T$ is a maximal simplex in $\operatorname{lk}x$ invisible from $v$.

For a four-face, $B$ is a $2$-cycle of triangles. Lemma~\ref{lem:local} applied in $\operatorname{lk}x$ shows that each $x$ extends at most one triangle of $B$, and that $v$ has degree three in this link. The same three triangles through the edge $vx$ give $x$ degree three in $\operatorname{lk}v$. Hence
\[
 4\le |B|\le n_3(\operatorname{lk}v)\le4 ,
\]
where $n_3$ counts cubic vertices. The first bound holds because the three edges of any triangle of $B$ must be paired with three distinct triangles. In the equality case, call the original triangle $pqr$. The other triangles have the form $pqs,prt,qru$. Pairing the edges $ps$ and $qs$ forces $s=t=u$. Thus $B$ is a tetrahedral boundary, and so is $\operatorname{lk}v$ by Lemma~\ref{lem:local}.

Let $A'=\{a_1,\dots,a_4\}=N(v)$ and $B'=\{b_1,\dots,b_4\}=V(B)$. These sets are disjoint since $B$ avoids $N[v]$. Equality in the count gives a bijection between the four triangles of $B$ and the four vertices of $A'$. Index them so that the simplex across $B'\setminus\{b_i\}$ is
\[
 \{a_i\}\cup(B'\setminus\{b_i\})\qquad(1\le i\le4).
\]
These are the four transversals that contain exactly one $a$, and each is a facet. The link of $a_i$ has seven vertices by Lemma~\ref{lem:local}. The seven neighbors of $a_i$ are therefore $v$, the three other vertices of $A'$ and the three vertices of $B'\setminus\{b_i\}$.

For the facets containing three $a$'s, fix $i$ and put $T_i=A'\setminus\{a_i\}$. Since $\operatorname{lk}v$ is a tetrahedral boundary, $\{v\}\cup T_i$ is a facet. Across $T_i$, the fourth vertex must be adjacent to all three vertices of $T_i$. The only candidates are $a_i$ and $b_i$, since these are the common neighbors other than $v$. If it were $a_i$, then $A'$ and the four facets through $v$ would form the boundary of a $4$-simplex. Each of its ridges would already belong to two of these facets. Ridge connectivity would then exclude the facets containing vertices of $B'$. Consequently the other four facets are
\[
 (A'\setminus\{a_i\})\cup\{b_i\}\qquad(1\le i\le4).
\]

This gives all eight choices of one vertex from each pair $(a_i,b_i)$ with an odd number of $a_i$. Choose a transversal uniformly. The selections from the four pairs are pairwise independent, since each assignment on two pairs has two completions. More generally a uniform fixed-parity $h$-bit vector is $(h-1)$-wise independent \cite[Corollary~6.14 and Example~6.16]{odonnell}. Let $\alpha_i,\beta_i$ be the signed columns on the walls labeled $a_i,b_i$. The vertex word $F$ is the fixed word of the four-face plus the four selected columns. Pairwise independence gives
\[
 R^2>\E\norm F^2=\norm{\E F}^2+\tfrac14\sum_{i=1}^4\norm{\alpha_i-\beta_i}^2>R^2 .
\]
The first inequality is \eqref{eq:corners}. On distinct axes $\alpha_i-\beta_i$ is an even word on two columns, so its norm exceeds $R$ by Proposition~\ref{prop:kernel}. On the same axis the two walls have opposite signs, giving norm $2>R$.

Hence the dual boundary of every four-face has no invisible maximal simplex. The link of a dual vertex is the dual boundary of the corresponding facet. The boundary argument in the first paragraph now gives the same assertion in every dimension above four by induction on face dimension.

Start with a pair of opposite facets of $K$, which exist by central symmetry. In a face of dimension at least four having both walls as facets, their dual vertices are nonadjacent. A maximal simplex through one must meet the closed neighborhood of the other. This gives a common neighbor, whose facet meets each wall in a facet of its own, by simplicity. Descending in this way gives a three-face with two facets on the same axis, contrary to Lemma~\ref{lem:axis}.
\end{proof}

\section{Triple supports in dimension three}
For $r=3$ each vertex saturates three axes. Counting edges of the dual identifies the icosahedron. Its signed faces give the five-cycle used in the Gram matrix.
Suppose $r=3$. An axis is active if it is saturated at a vertex.
\begin{prop}[Triple supports]\label{prop:five}
After column sign switches and relabeling, the active axes are $\{\infty,0,1,2,3,4\}$ and the ten vertex supports are
\begin{equation}\label{eq:faces}
 \{\infty,i,i+1\},\qquad (\mathbb Z/5)\setminus\{i,i+2\}\qquad (i\in\mathbb Z/5).
\end{equation}
Within a support, pairwise sign products are positive on pairs incident to $\infty$ and on cycle edges, and negative on cycle diagonals.
\end{prop}
\begin{proof}
Let $m\ge3$ be the number of active axes. The dual graph $G$ of $K$ triangulates the sphere and has $2m$ vertices. Each has degree at most five by Lemma~\ref{lem:axis}. By Lemma~\ref{lem:sign}, edges on the same axis pair have the same relative sign, so there are at most two per pair. Euler gives
\[
 6(m-1)=|E(G)|\le\min\{5m,m(m-1)\}.
\]
Thus $m=6$ and every degree is five. The dual is the icosahedron \cite[Table~1]{bm}. Opposite walls share no neighbor by the same relative-sign condition, so they are antipodes.

Write $i^\pm$ for the walls on axis $i$. Switch column signs to label one pole $\infty^+$ and its cyclic neighbors $u_i=i^+$ ($i\in\mathbb Z/5$). Let $l_i$ be the third vertex of the triangle across the edge $u_i u_{i+1}$ from $\infty^+$. The $l_i$ form the lower pentagon around $\infty^-$. The edges of the two pentagons and their belt give the neighbor sets
\[
\begin{aligned}
 N(u_i)&=\{\infty^+,u_{i-1},u_{i+1},l_{i-1},l_i\},\\
 N(l_{i+2})&=\{\infty^-,l_{i+1},l_{i+3},u_{i+2},u_{i+3}\}.
\end{aligned}
\]
Thus $u_i$ and $l_{i+2}$ are nonadjacent and have no common neighbor, so they are antipodes. In particular $l_i=(i+3)^-$. Up to antipodes the faces are the five polar triangles $(\infty^+,u_i,u_{i+1})$ and the five belt triangles $(u_i,u_{i+1},l_i)$. Their labels give \eqref{eq:faces} and the stated sign products.
\end{proof}

\section{A rank bound for the five-cycle}
Fix the block signs, put $z=\sum_b\theta_b E_b$, and let $q=\norm z^2$. For a vertex word $f$, average the bounds on $\norm{f\pm z}^2$ in \eqref{eq:corners} to obtain
\begin{equation}\label{eq:q}
  \norm f^2<1-q,\qquad 0\le q<1 .
\end{equation}
Let $U$ consist of the six active columns, $W$ of the remaining $\ell=d-3$ columns, and $f_T\in\{0,\pm1\}^6$ be a vertex sign vector on $T$. Let $\Sigma$ be the zero-diagonal matrix of pairwise sign products, $H=I+\Sigma$, $B=U\T U$, and $y_j=UHe_j$. Each row of $\Sigma$ has five entries of magnitude one, so $(\Sigma^2)_{ii}=5$. By symmetry of the matrix and rotation of the five-cycle, it suffices to check three off-diagonal entries. Summing over the four remaining indices in the order $(\infty,0,1,2,3,4)$ gives
\[
 (\Sigma^2)_{\infty,0}=1-1-1+1,\qquad
 (\Sigma^2)_{0,1}=1-1+1-1,\qquad
 (\Sigma^2)_{0,2}=1+1-1-1 .
\]
All three sums vanish, so $\Sigma^2=5I$. Set
\[
 p_T=1-q-\norm{Uf_T}^2,\qquad \delta_j=\norm{y_j}^2+q-2,\qquad \bar\delta=\tfrac16\sum_{j=1}^6\delta_j .
\]
Put $A_0=(W,z)$, $C=A_0\T A_0$, $X=A_0\T UH$ and $X_j=A_0\T y_j$. Column $j$ of $H$ assigns signs to the columns of $U$. Thus every trailing sign vector $\tau$ satisfies
\[
 L<\min_\pm\norm{A_0\tau\pm y_j}^2
 =\tau\T C\tau+\norm{y_j}^2-2|\ip{X_j}{\tau}|.
\]
Taking fair $\tau$ gives $\delta_j>2\E|\ip{X_j}{\tau}|\ge\sqrt2\norm{X_j}$ by the real Khintchine inequality \cite{sz}, since $\tr C=\ell+q$. Dropping the absolute-value term and averaging over $j$ gives $\tau\T C\tau>\tr C-\bar\delta$. Hence \eqref{eq:moment} gives
\begin{equation}\label{eq:tail}
 C\succeq\diag C-\tfrac12\bar\delta I,\qquad
 \diag C=I-(1-q)ee\T ,
\end{equation}
where $e$ is the last trailing coordinate vector.

\begin{lemma}[The five-cycle Gram matrix]\label{lem:aff}
Put $D=HBH$ and $u_T=(\Sigma-I)f_T$. Each $u_T$ is a sign vector, and $D$ is a convex combination of the following leading blocks.
\begin{equation}\label{eq:aff}
\begin{array}{lccc}
 \text{weight} & q & p_T/4 & \delta_j/4 \\
 \text{leading block} & I-\Sigma &4I-2u_Tu_T\T&2I+4e_je_j\T .
\end{array}
\end{equation}
\end{lemma}
\begin{proof}
Put $K_0=\Sigma-I$ and use the signed triples of Proposition~\ref{prop:five}. For $a=e_j$ or $f_T$, the entries of $K_0a$ are odd and $\norm{K_0a}^2=6\norm a^2-2a\T\Sigma a=6$, since $f_T\T\Sigma f_T=6$. Thus $K_0a$ is a sign vector. Choose its global sign so that the $\infty$ coordinate is positive. Write $N_-(K_0a)$ for the resulting set of negative coordinates. Multiplying by $K_0$ gives the four rows below, with $i$ taken modulo five.
\[
\begin{array}{lc}
 a & N_-(K_0a)\\ \hline
 e_\infty & \mathbb Z/5\\
 e_i & \{i,i+2,i+3\}\\
 e_\infty+e_i+e_{i+1} & \{i+3\}\\
 e_i+e_{i+1}-e_{i+3} & \{i+2,i+3,i+4\}
\end{array}
\]
For instance, in the order $(\infty,0,1,2,3,4)$ the last row at $i=0$ gives $(1,1,1,-1,-1,-1)$. The three-element negative sets in the last row have cycle edges as complements. Those in the $e_i$ row have cycle diagonals as complements. Together the rows give every odd subset of $\mathbb Z/5$ once. Hence the sixteen vectors represent the full six-coordinate odd-parity class up to negation.

The second and fourth moments on this class agree with those of independent fair signs by the fixed-parity observation in the proof of Proposition~\ref{prop:high}. These moments are invariant under global negation, so one may average over the sixteen representatives. The identity $HK_0=4I$ gives $(K_0a)\T D(K_0a)=16\norm{Ua}^2$. For $a=e_j$ this is $16$, and for $a=f_T$ it is $16(1-q-p_T)$. Substituting in \eqref{eq:moment} with $t=16$ yields
\[
 2(D-\diag D)+(\tr D-16)I=-\sum_T(q+p_T)u_Tu_T\T .
\]
The diagonal gives $\sum_T p_T+\sum_j\delta_j=4(1-q)$, using $D_{jj}=2-q+\delta_j$. By the second-moment identity, $16I$ is the sum of the sixteen outer products. The six vectors $K_0e_j$ contribute $K_0^2=6I-2\Sigma$, leaving
\[
 \sum_T u_Tu_T\T=16I-K_0^2=10I+2\Sigma .
\]
Substitution in the off-diagonal entries gives \eqref{eq:aff}. The weights sum to one by the diagonal identity, and $p_T>0$ by \eqref{eq:q}.

\end{proof}

\begin{lemma}[The five-cycle rank bound]\label{lem:spec}
Under the preceding hypotheses, $Y=(UH,2A_0)$ satisfies $\dim\ker Y\le3$.
\end{lemma}
\begin{proof}
Complete the $p_T$- and $\delta_j$-blocks of \eqref{eq:aff} to symmetric matrices $M_T,J_j$ with trailing blocks $4I,2I$ and lower-left blocks $0,8X_je_j\T/\delta_j$, respectively. With $Q=q(I-\Sigma)\oplus(4qI-4(1-q)ee\T)$,
\[
 Y\T Y\succeq Q+\tfrac14\Big(\sum_T p_T M_T+\sum_j\delta_j J_j\Big).
\]
The leading and mixed blocks agree. The trailing block on the right is $(4-3\bar\delta)I-4(1-q)ee\T\preceq4C$ by \eqref{eq:tail}.

Write $s_4$ for the sum of the four smallest eigenvalues. Then $s_4(M_T)=-8+3\cdot4=4$. On the orthogonal complement of a plane $J_j=2I$. On that plane $J_j+2I$ has diagonal $8,4$ and determinant $32-64\norm{X_j}^2/\delta_j^2>0$. Thus $\lambda_{\min}(J_j)>-2$ and the larger eigenvalue on the plane is at least six, giving $s_4(J_j)>4$.

By \eqref{eq:aff}, $0\preceq D\preceq q(I-\Sigma)+6(1-q)I$, so $q\le6/(5+\sqrt5)<4/(7-\sqrt5)$. Thus $8q-4\le q(1+\sqrt5)\le4q$. Since $\Sigma^2=5I$ and $\tr\Sigma=0$, the four lowest eigenvalues of $Q$ are $8q-4$ and three copies of $q(1-\sqrt5)$. For a real symmetric matrix $A$, Ky Fan's minimum principle \cite{fan} gives
\[
 s_4(A)=\min_{\substack{P=P^\mathsf T=P^2\\ \operatorname{rank}P=4}}\tr(PA).
\]
Choose a minimizing projection $P$ for $Y\T Y$. Taking traces after multiplication by $P$ preserves the block inequality. For each matrix $B$ on the right we have $\tr(PB)\ge s_4(B)$ for this same $P$. Summing with the nonnegative weights gives
\[
\begin{aligned}
 s_4(Y\T Y)
 &\ge s_4(Q)+\tfrac14\Big(\sum_T p_T s_4(M_T)+\sum_j\delta_j s_4(J_j)\Big)\\
 &> s_4(Q)+4(1-q)=(7-3\sqrt5)q\ge0 .
\end{aligned}
\]
Here $\sum_T p_T+\sum_j\delta_j=4(1-q)$ by \eqref{eq:aff}. The strict inequality follows from $\delta_j>0$ and $s_4(J_j)>4$, including when $q=0$. Thus $\dim\ker Y\le3$.
\end{proof}

\begin{proof}[Proof of Theorem~\ref{th:main}]
For $d\ge3$, Proposition~\ref{prop:kernel} gives a positive odd integer $r$. If $r=1$, a single unit column violates \eqref{eq:corners}. Proposition~\ref{prop:high} excludes $r\ge5$. For $r=3$, Proposition~\ref{prop:five} gives the five-cycle. Its matrix $Y$ has $d$ rows and $d+4$ columns, so $\dim\ker Y\ge4$, contrary to Lemma~\ref{lem:spec}.

For $d=1$, the even number of unit terms can be signed to cancel in pairs. For $d=2$, embed the odd family in the horizontal coordinate plane of $\R^3$ and append the unit normal $e_3$. The three-dimensional bound gives $D(V)+1=D((V,e_3))\le2$ by orthogonality.
\end{proof}

\section*{Declaration of generative AI and AI-assisted technologies in the manuscript preparation process}
OpenAI Codex assisted with the English translation. The authors take full responsibility for the manuscript.

\end{document}

%% file: metadata.tex
\newcommand{\PaperTitle}{A sharp parity bound for signed sums of Euclidean unit vectors}
\newcommand{\LiName}{Muwen Li}
\newcommand{\ZhaoName}{Chunyin Zhao}
\newcommand{\LiAff}{School of Mathematical Sciences, University of Science and Technology of China}
\newcommand{\ZhaoAff}{Beijing Technology and Business University}
\newcommand{\LiPostal}{No.~96 Jinzhai Road, Hefei, Anhui 230026, P.R. China}
\newcommand{\ZhaoPostal}{Liangxiang Higher Education Park, Fangshan District, Beijing 102488, P.R. China}